\documentclass[a4paper,10pt]{article}

\usepackage[utf8x]{inputenc}
\usepackage[english]{babel}
\usepackage[T1]{fontenc}
\usepackage{amsmath,amscd}
\usepackage{amsfonts}
\usepackage{amssymb}
\usepackage{amsthm}
\usepackage{mathrsfs}
\usepackage{graphicx}
\usepackage[all,cmtip]{xy}
\usepackage{textcomp}
\usepackage{qtree}
\usepackage{url}
\usepackage{hyperref}
\hypersetup{colorlinks=true,linkcolor=blue}

\newtheorem{teo}{Theorem}[section]
\newtheorem{lemma}{Lemma}[section]
\newtheorem{rem}{Remark}[section]
\newtheorem{prop}{Proposition}[section]
\newtheorem{cor}{Corollary}[section]

\newtheorem{defin}{Definition}[section]

\newcommand{\Z}{{\mathbb{Z}}}
\newcommand{\C}{{\mathbb{C}}}
\newcommand{\R}{{\mathbb{R}}}
\newcommand{\Q}{{\mathbb{Q}}}
\newcommand{\N}{{\mathbb{N}}}

\newcommand{\G}{{\mathcal{G}}}
\newcommand{\Proj}{{\mathbb{P}}}
\newcommand{\A}{{\mathcal{A}}}
\newcommand{\K}{{\mathbb{K}}}
\newcommand{\Sc}{{\mathcal{S}}}
\newcommand{\Cc}{{\mathcal{C}}}
\newcommand{\F}{{\mathcal{F}}}
\newcommand{\T}{{\mathcal{T}}}
\newcommand{\U}{{\mathcal{U}}}
\newcommand{\Hc}{{\mathcal{H}}}

\newcommand{\maxAnmenouno}{{\Cc_{A_{n-1}}}}
\newcommand{\maxBn}{{\Cc_{B_n}}}
\newcommand{\Lev}{{\text{Lev}}}
\newcommand{\rk}{{\text{rk}}}
\newcommand{\Aut}{{\text{Aut}}}
\newcommand{\out}{{\text{out}}}
\newcommand{\ptilde}{{\widetilde{p}}}
\newcommand{\tr}{{\text{tr}}}
\newcommand{\ptildetilde}{{\widetilde{\ptilde}}}
\newcommand{\ptildetildeAtras}{{\ptildetilde_A\left[1\right]}}
\newcommand{\wout}{{\text{wout}}}
\newcommand{\Qtilde}{{\widetilde{Q}}}
\newcommand{\Qtildetilde}{{\widetilde{\Qtilde}}}

\newcommand{\Qtildetildetras}{{\Qtildetilde_B\left[1\right]}}
\newcommand{\maxDn}{{\Cc_{D_n}}}
\newcommand{\supp}{{\text{supp}}}
\newcommand{\cont}{{\text{cont}}}

\title{Poincar\'e series for maximal De Concini-Procesi models of root arrangements}
\author{Giovanni Gaiffi and Matteo Serventi}
\date{\today}

\begin{document}

\maketitle




\section{Introduction}
In this paper we focus on  maximal complex De Concini-Procesi  models
associated to root arrangements of types A, B, C, D and we compute inductive formulas for their Poincar\'e series.

In \cite{DCP2}, \cite{DCP1},   De Concini and Procesi constructed  {\em wonderful  models}  for the   complement of a
subspace arrangement in a  vector space.  In general, given a subspace arrangement, there are several De
Concini-Procesi
models associated to it,   depending on   
  distinct sets 
of initial combinatorial data (``building sets'', see Section \ref{sec:building}). 

The  interest in these  varieties was at first motivated by an  approach to Drinfeld construction of special
solutions for Khniznik-Zamolodchikov equation (see \cite{drinfeld}). Then real and complex  De Concini-Procesi models
turned out to play a central role in several   fields of mathematical research:  subspace (and toric)  arrangements,  toric varieties,
moduli spaces of curves,     configuration spaces, box splines, index theory  and discrete  geometry    (see for instance   \cite{DCP4}, \cite{DCP3}, 
\cite{etihenkamrai}, \cite{feichtner}, \cite{feichtnerkozlov},  \cite{gaiffi2},   \cite{postnikov},
\cite{postnikoreinewilli} and  \cite{zelevinski}).

 Among the building sets  associated to a given subspace arrangement there are  always a minimal one and a maximal one
with respect to inclusion: as a consequence there are always a  minimal and a maximal De Concini-Procesi model. 
Several examples of minimal   models (associated to the minimal   building set of {\em irreducible subspaces})  have been studied
in detail. More recently, the relevance of real and complex maximal models was pointed out (maximal models appear for
instance in   \cite{LTV},  \cite{SzenesVergne},  and in the context  of   toric varieties, see  \cite{HenPisa} for
further references). 

The case of root arrangements is particularly interesting. 
Let us consider for instance the arrangements of type \(A_n\). Our purpose is to compute the series
\[
 \phi_A\left(q,t\right)=t+\sum_{n\geq 2} P_\maxAnmenouno\left(q\right)\frac{t^n}{n!} \ \in\Q[q][[t]]
\]
Here, for every $n\geq 2$,
$$P_\maxAnmenouno\left(q\right)=\sum_{}^{}\left(\rk\left(H^{2i}\left(Y_\maxAnmenouno,\Z\right)\right)\right)q^{i}
$$
is the Poincar\'e polynomial of the maximal complex De Concini-Procesi model $Y_\maxAnmenouno$ and the variable $q$ has degree two (in odd degree, the integer cohomology of De Concini-Procesi models is 0  - see \cite{DCP1}).

By  carefully  counting  the elements  of a   basis for the integer cohomology  which was first described by Yuzvinski (see
\cite{YuzBasi} and also \cite{GaiffiBlowups}),   we find 
 an inductively defined  series in infinite variables \(g_0, g_1, g_2, \ldots g_n,\ldots\) with the following
property:  when we replace  $g_0$ with $t$ and, for every $i,r \ \geq 1$, $g_i^r$ with
$\frac{q^r-q}{q-1}t^r$, we obtain the series $\phi_A\left(q,t\right)$ (Theorem \ref{serie da ptilde}). 

Some explicit computations (using  the Computer Algebra system Axiom)  show that our method is effective (see Section
\ref{sec:examples}).
The same technique can be extended to the case  \(B_n\) (see Theorem \ref{serie da Q}), which, from the point of view of
subspaces and models, is equal to \(C_n\), and also to the \(D_n\) case (see Section
\ref{subsection:Dn}).

In Section  \ref{sec:final} 
we  show  that the  series in infinite variables  computed in the preceding sections encode  more general results.  For instance,  they  allow us to obtain  the Poincar\'e series of the   families  of   De Concini-Procesi  models    whose building sets are  the maximal building sets \(\Cc_{A_n}, \Cc_{B_n},\Cc_{D_n}\) tensored  by \(\C^h\):  it is sufficient to perform  different substitutions of  the variables \(g_0, g_1, \ldots g_n,\ldots\)
We observe that in the \(A_n\) case,  the complements of these tensored arrangements  are classical generalizations of the pure braid space  (see  \cite{cohen1976hil} and  \cite{sundaramwelker}).


\section{Basic Concepts}

\subsection{Some combinatorics of subspace arrangements}
\label{sec:building}
Let $V$ be a finite dimensional vector space over an infinite field $\K$ and denote by $V^*$ its dual.
Let now $\G$ be a finite set of subspaces of $V^*$ and denote by $\Cc_\G$ its closure under the sum.

\begin{defin}
 Given a subspace $U\in\Cc_\G$, a \textbf{decomposition of} $U$ in $\mathbf{\Cc_\G}$ is a collection
$\{U_1,\cdots,U_k\}$ ($k>1$) of non zero subspaces in $\Cc_\G$ such that
\begin{enumerate}
 \item $U=U_1\oplus\cdots\oplus U_k$
 \item for every subspace $A\subset U$, $A\in\Cc_\G$, we have $A\cap U_1,\cdots,A\cap U_k \in \Cc_\G$ and
$A=\left(A\cap U_1\right)\oplus\cdots\oplus \left(A\cap U_k\right)$.
\end{enumerate}
\end{defin}
\begin{defin}
 A subspace $F\in\Cc_\G$ which does not admit a decomposition  is called \textbf{irreducible} and the set of
irreducible subspaces is denoted by $\mathbf{\F_\G}$.
\end{defin}
One can prove (see  \cite{DCP1}) that 
every subspace $U\in\Cc_\G$ has a unique decomposition into irreducible subspaces.

\begin{defin}\label{building}
 A collection $\G$ of subspaces of $V^*$ is called \textbf{building} if every element $C\in\Cc_\G$ is the direct sum
$G_1\oplus\cdots\oplus G_k$ of the set of maximal elements $G_1,\cdots,G_k$ of $\G$ contained in $C$.
\end{defin}
\begin{rem}(see \cite{DCP1})
 \begin{itemize}
  \item The set of irreducible subspaces  of a given family of subspaces of $V^*$ is building.
  \item A set  of subspaces of $V^*$ which is closed under the sum is building.
 \end{itemize}
\end{rem}
Given a family $\G$ of subspaces of $V^*$ there are different sets \({\mathcal B}\) of subspaces of $V^*$ such that
\(\Cc_{{\mathcal B}}=\Cc_\G\); if we order by inclusion the collection of such sets, it turns out that the minimal
element is
$\F_\G$ and the maximal one is $\Cc_\G$.
\begin{defin}(see \cite{DCP1}) \label{Gnested}
 Let $\G$ be a building set of subspaces of $V^*$. A subset $\Sc\subset \G$ is called $\mathbf{\G}$\textbf{-nested} if and only if 
 for every  subset $\{A_1,\cdots,A_k\}$ (\(k\geq 2\)) of pairwise non comparable elements of $\Sc$  the subspace   $A=A_1+\cdots +
A_k$ does not belong to  $\G$.

\end{defin}
\begin{rem}\label{Gnested sets nel caso massimale}
 If  $\Cc$ is a building family of subspaces closed under the sum, then   the subspaces of 
a $\Cc$-nested set  must be totally ordered (with respect to inclusion).
\end{rem}

\subsection{Wonderful models}

Let now $V$ be a finite dimensional complex vector space and denote by $V^*$ its dual.\\
Let us consider a finite subspace arrangement $\G$ in $V^*$ and, for every $A\in\G$, let us denote by  $A^\perp$ its
annihilator in $V$.\\
For every $A\in\G$ we have a rational map
$$\pi_A:V\longrightarrow V/A^\perp \longrightarrow \Proj \left( V/A^\perp \right)$$
which is regular on $V-A^\perp$.\\
We then consider the embedding 
$$\phi_\G:\A_\G\longrightarrow V\times\prod_{A\in\G}\Proj\left( V/A^\perp \right)$$
given by inclusion on the first component and the maps \(\pi_A\) on the other components. The De Concini-Procesi model
\(Y_{\G}\) associated to  $\G$ is the closure of $\phi_\G \left(\A_\G \right)$ in
$V\times\prod_{A\in\G}\Proj\left(V/A^\perp\right)$.
\medskip

These  {\em wonderful models}  are  particularly interesting when the arrangement $\G$ is building: they turn out to be 
smooth varieties and the complement of \(\A_\G\) in \(Y_\G\) is a divisor with normal crossings, described in terms of
\(\G\)-nested sets.  Moreover, their integer cohomology rings
are  torsion free    (see \cite{DCP1}). In \cite{YuzBasi}  Yuzvinski explicitly described  $\Z$-bases of these rings
(see also  \cite{GaiffiBlowups}): we  briefly  recall these  results concerning cohomology.

Let $\G$ be a building set of subspaces of $V^*$. Let $\Hc\subset
\G$, let $B\in\G$ such that $A\subsetneq B$ for each $A\in\Hc$ and define
$$d_{\Hc,B}:=\dim B - \dim \left(\sum_{A\in\Hc} A\right).$$
With these notations, in the polynomial ring $\Z[c_A]_{A\in\G}$, we put 
$$P_{\Hc,B}:=\prod_{A\in\Hc}c_A\left(\sum_{C\supset B}c_C\right)^{d_{\Hc,B}}$$
and we denote by $I$ the ideal generated by these polynomials as $\Hc$ and $B$ vary.
\begin{teo}(see \cite{DCP1}).\\
There is a surjective ring homomorphism
$$\phi \: : \: \Z[c_A]_{A\in\G}\longrightarrow H^*(Y_\G,\Z)$$
with  kernel is $I$ and such that  \(\phi(c_A)\in H^2(Y_\G,\Z)\).  
\end{teo}
\begin{defin}
 Let $\G$ be a building set of subspaces of $V^*$. A function
$$f:\G\longrightarrow \N$$
is $\mathbf{\G}$\textbf{-admissible} (or simply \textbf{admissible}) if $f=0$ or, if $f\neq 0$, $\supp(f)$ is
$\G$-nested and for all $A\in\supp(f)$ one has
$$f(A)< d_{\supp(f)_A,A}$$
where $\supp(f)_A:=\{C\in\supp(f):C\subsetneq A\}$.
\end{defin}
\begin{defin}
 A monomial $m_f=\prod_{A\in\G}c_A^{f(A)}\in\Z[c_A]_{A\in\G}$ is \textbf{admissible} if $f$ is admissible.
\end{defin}
\begin{teo}\label{base coomologia}(see \cite{YuzBasi} and also \cite{GaiffiBlowups})\\
 The set $\mathcal{B}_\G$ of all admissible monomials corresponds to  a $\Z$-basis of $H^*(Y_\G,\Z)$.
\end{teo}

\section{Maximal Models of Reflection Arrangements of Classical Type}

We  now focus on  the cohomology rings of maximal De Concini-Procesi models for  
root  arrangements  of type $A_n$, $B_n$ (\(C_n\))  and $D_n$.

\subsection{Type $A_{n-1}$}

Let $W$ be a complex vector space of dimension $n$ and consider the arrangement given by hyperplanes
$H_{ij}:=\{z_i-z_j=0\}$ where $z_i$ ($i=1,\cdots,n$) are the coordinates. The intersection of these
hyperplanes is the line $N=\{z_1=\cdots=z_n\}$:   we consider the quotient $V=W/N$ and the arrangement provided by the
images of the hyperplanes $H_{ij}$ via the quotient map $W\stackrel{\pi}{\longrightarrow}V$.\\
We can choose     linear functionals  $f_{ij}$  in $V^*$ such that  the zeroes of $f_{ij}$  form the hyperplane $\pi\left(H_{ij}\right)$
and the set \(\{f_{ij}\}\) is    a root system of type $A_{n-1}$.\\
In $V^*$ we consider the subspace arrangement $\A_{A_{n-1}}$ given by the lines $<f_{ij}>$ and denote for brevity by  $\maxAnmenouno$ its
closure under the sum and by $\F_{A_{n-1}}$ the set of irreducible subspaces  in $\maxAnmenouno$.\\
Our purpose is to compute the series
\begin{equation}\label{serie poincare' An}
 \phi_A\left(q,t\right)=t+\sum_{n\geq 2} P_\maxAnmenouno\left(q\right)\frac{t^n}{n!} \ \in\Q[q][[t]]
\end{equation}
where, for every $n\geq 2$,
$$P_\maxAnmenouno\left(q\right)=\sum_{}^{}\left(\rk\left(H^{2i}\left(Y_\maxAnmenouno,\Z\right)\right)\right)q^{i}
$$
is the Poincar\'e polynomial of $Y_\maxAnmenouno$ (the variable $q$ has  degree 2).

In \cite{YuzBasi}, Yuzvinsky noticed  there is a bijective correspondence (actually an isomorphism of
partially ordered sets) between the elements of $\F_{A_{n-1}}$ and the subsets of $\{1,\cdots,n\}$ of cardinality at
least two:  this correspondence identifies the subset $\{i_1,\cdots,i_k\}$   with 
$<f_{i_1i_2},\cdots,f_{i_{k-1}i_k}>$.
Since every subspace in $\maxAnmenouno$ has a unique decomposition into irreducible subspaces, we can identify
elements of $\maxAnmenouno$
with families of disjoint subsets of cardinality at least two of $\{1,\cdots,n\}$. Furthermore, given two such
collections
$X=\{X_1,\cdots,X_k\}$ and $Y=\{Y_1,\cdots,Y_r\}$ we say that $Y$ is included in $X$ (and write $Y\subset X$) if for every
$i\in\{1,\cdots,r\}$ there exists $j\in\{1,\cdots,k\}$ such that $Y_i\subset X_j$; our identification thus becomes an
isomorphism of partially ordered sets (we order $\maxAnmenouno$ by inclusion).

A  $\maxAnmenouno$-nested set is a subset  of   $\maxAnmenouno$ strictly ordered by inclusion.  Now we will see  how to associate a graph (actually a forest of levelled oriented rooted trees) to a
$\maxAnmenouno$-nested set.
Let us  first recall by  an example the $\F_{A_{n-1}}$ case  (see 
\cite{YuzBasi}).   Let us take, as
$\F_{A_{8}}$-nested set, the collection $\Sc:=\{(1,2,3,4,5),(6,7,8,9),(1,4,5),(6,7)\}$; we associate to $\Sc$ the
following graph:
$$
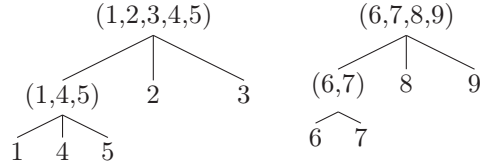
\qtreecenterfalse
\Tree [ .(1,2,3,4,5) [ 1 4 5 ].(1,4,5) 2 3 ]
\hskip 0.3in
\Tree [ .(6,7,8,9) [ 6 7 ].(6,7) 8 9 ]
$$
where the edges are directed from the top to the bottom.

Let us now consider  the $\Cc_{A_{15}}$-nested set
$$\Sc:=\{(1,2,3,4,5)(8,10,12,13,14,16),(1,2,4,5)(8,10,12),(1,2),(10,12)\}.$$   We associate to it the {\em levelled} forest  $\Gamma\left(\Sc\right)$: 
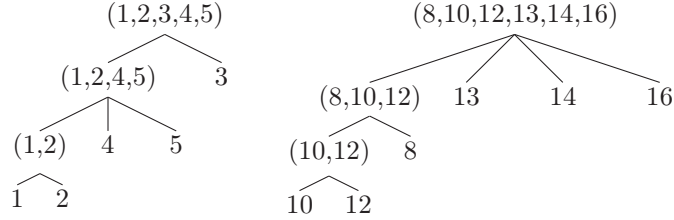
$$
\qtreecenterfalse
\Tree [ .(1,2,3,4,5) [ [ 1 2 ].(1,2) 4 5 ].(1,2,4,5) 3 ]
\hskip 0.3in
\Tree [ .(8,10,12,13,14,16) [ [ 10 12 ].(10,12) 8 ].(8,10,12) 13 14 16 ]
$$
where, again, the orientation is from the top to the bottom, and the elements of the nested set can be read  ``level by level''  from the vertices which are not leaves (we call level 1 the level which contains the roots, 
and level $k+1$ the one which  contains the  vertices which are $k$ steps away from a  root).



Let  now   \(\Sc\)  be  a $\maxAnmenouno$-nested set and  let us denote by \(B\) (resp. \(A\)) the element of \(\Sc\) determined by the vertices (not leaves) at level $k$ (resp. $k+1$). Then   $A$ is the maximal element of $\Sc$ strictly contained
in $B$.\\
Hence if $B$ is given by the family $\{B_1,\cdots,B_k\}$, $A$ by $\{A_1,\cdots,A_r\}$ and, for every
$i\in\{1,\cdots,k\}$, we set $I_{B_i}:=\{j\in\{1,\cdots,r\}:A_j\subset B_i\}$, we have
$$\dim B_i - \sum_{j\in I_{B_i}}\dim A_j= |\out(v_{B_i})|-1$$
where $v_{B_i}$ is the vertex of $\Gamma\left(\Sc\right)$ associated to  $B_i$ and $\out(v_{B_i})$ is the set of outgoing
edges from $v_{B_i}$. Then we obtain: 
\begin{equation}\label{bound coomologia in termini di grafi}
 d_{\{A\},B}=\sum_{v\in \Lev(k)}\left(|\out(v)|-1\right)
\end{equation}
where $\Lev(k)$ is the set of vertices (not leaves) of $\Gamma\left(\Sc\right)$ belonging to level $k$.
\begin{defin}\label{foreste ammissibili}
 A levelled forest $\Gamma$ is \textbf{admissible} if, for any of its level $k$, one has
$$\sum_{v\in\Lev(k)}\left(|\out(v)|-1\right)=\left(\sum_{v\in\Lev(k)}|\out(v)|\right)-|\Lev(k)|\geq 2.$$ 
\end{defin}

\begin{defin}
 Let $\Gamma$ be a levelled admissible forest on $n\geq 2$ leaves. We denote by \(\cont(\Gamma)\)  the  contribution given
 to the series (\ref{serie poincare' An})  by all  the monomials \(m_f\) of the basis  such that \(\supp \ f\) is a
nested set whose  graph is (up to a relabelling of the leaves) isomorphic to \(\Gamma\).
\end{defin}
\begin{prop}\label{contributo foresta ammissibile}
 Let $\Gamma$ be a levelled admissible forest on $n\geq 2$ leaves.   Then we have 
$$\cont(\Gamma)=\frac{n!}{|\Aut\left(\Gamma\right)|}C_\Gamma(q)\frac{t^n}{n!}$$
where
$$C_\Gamma(q)=\prod_{k \text{ level}}\frac{q^{\sum_{v\in\Lev(k)}\left(|\out(v)|-1\right)}-q}{q-1}$$
and $\Aut\left(\Gamma\right)$ is the group of automorphisms of $\Gamma$.
\end{prop}
\proof
We notice that 
there are $\frac{n!}{|\Aut\left(\Gamma\right)|}$ different $\maxAnmenouno$-nested sets whose associated graph is
$\Gamma$.\\
The thesis follows by observing that, if $\Gamma=\Gamma\left(\Sc\right)$ where \(\Sc\) is a $\maxAnmenouno$-nested set,
the  contribution to the Poincar\'e polynomial $P_\maxAnmenouno(q)$
of the monomials \(m_f\)  such that \(\supp \ f= \Sc\) is $C_\Gamma(q)$.\\
\qed 

Our idea to compute the series (\ref{serie poincare' An}) is to consider all levelled forests (not necessarily
admissible) on at least two leaves and associate, to each of them, a monomial that encodes data we are interested in:
number of levels and, for each level \(k\), the number \(\sum_{v\in \Lev(k)}\left(|\out(v)|-1\right)\). We will put
together these monomials in a series, which will be calculated inductively, and from which one can obtain the series
(\ref{serie poincare' An}).
\begin{defin}\label{coda banale}
 A  \textbf{trivial tail} of a (levelled) oriented forest  $\T$  is  given by a subtree $\T'$
which stems from  a  vertex \(v\)  of $\T$  with $|\out(v)|=1$ and has a single leaf. 
\end{defin}
\begin{defin}
 Two (levelled) oriented forests  $\T_1$ and $\T_2$ are \textbf{equivalent} if they differ only for trivial tails.
\end{defin}
\begin{defin}
 Given an equivalence class of (levelled) rooted oriented forests  modulo trivial tails we call \textbf{minimal representative}
the tree in this class with no trivial tails.
\end{defin}

Let us now define the following series
$$\ptilde_A(g_0,g_1,g_2,g_3,\cdots):=g_0+\sum_\Gamma\frac{g_0^{\tr\left(\Gamma\right)}}{|\Aut\left(\Gamma\right)|}
\prod_v g_{l(v)}^{|\out(v)|-1}
$$
where $\Gamma$ runs among minimal representatives of levelled oriented forests on at least two leaves,
$\tr\left(\Gamma\right)$ is the
number of trees of  $\Gamma$ (we are considering also the degenerate tree  given by a single leaf), $v$ varies
among the vertices (not leaves) of $\Gamma$ and $l(v)$ is the level of $v$.
\begin{defin}
 A monomial of $\ptilde_A$ is \textbf{bad} if there exist $1\leq i<j$ such that $g_j$ appears in the monomial but $g_i$
doesn't.
\end{defin}
We notice that bad monomials correspond to forests with a level \(k\)   such that \(\sum_{v\in \Lev(k)}\left(|\out(v)|-1\right)=0\).

\begin{defin}
 A monomial $m$ of $\ptilde_A$ has \textbf{valency} $k$ if $k=\max\{j\geq 1: g_j \text{ appears in } m\}$.
\end{defin}
\begin{prop}\label{grado e foglie nelle foreste}
 Given  a levelled oriented  forest $\Gamma$ on $n\geq 2$ leaves,   let $m_\T$ be the monomial of $\ptilde_A$ associated to $\T$.
Then the degree of $m_\Gamma$  is \(n\).
\end{prop}
\proof
One can restrict to trees and then proceed by   induction on the valency of $m_\Gamma$.
\qed

\begin{teo}\label{serie da ptilde}
 Removing bad monomials from $\ptilde_A$ and replacing $g_0$ with $t$ and, for every $i,r \ \geq 1$, $g_i^r$ with
$\frac{q^r-q}{q-1}t^r$ we obtain the Poincar\'e series (\ref{serie poincare' An}). 
\end{teo}
\proof
Let's start by observing that, if we remove bad monomials, we haven't removed all  monomials corresponding to non admissible forests; in fact we
still have the ones corresponding to forests  in which there is (at least) a level, say $k\geq 1$, such that
$$\sum_{v\in\Lev(k)}\left(|\out(v)|-1\right)=1.$$
Anyhow the contribution of  such monomials is  killed by our substitution; indeed they  have a variable whose exponent  is $1$ and
$\frac{q^r-q}{q-1}=0$ if $r=1$.\\
Let now $\Gamma$ be an admissible forest on $n\geq 2$ leaves; let  $m_\Gamma$ be the
monomial of $\Gamma$ in $\ptilde_A$ and  $k\geq 1$ be its valency. For all $1\leq j\leq k$ the exponent of the
variable $g_j$ in $m_\Gamma$ is
$$\sum_{v\in \Lev(j)}\left(|\out(v)|-1\right).$$
Our claim then follows from Proposition  \ref{contributo foresta ammissibile} and Proposition \ref{grado e foglie nelle foreste}.\\
\qed
\\
The problem is now reduced to the computation of  $\ptilde_A$. To this end we define
$$\ptildetilde_A:=1+\sum_\T \frac{1}{|\Aut\left(\T\right)|}\prod_v g^{|\out(v)|-1}_{l(v)}$$
where $\T$ runs among minimal representatives of levelled oriented trees on $n\geq 2$ leaves and $v$ among the vertices
(not leaves) of $\T$.
 One immediately checks that:
\begin{equation}
\label{ptilde in termini di ptildetilde}
\ptilde_A=e^{g_0\ptildetilde_A}-1,
\end{equation}
therefore,  all we need is a formula for  $\ptildetilde_A$.
\begin{teo}\label{formula induttiva per ptildetildeA teorema}
 The following recursive  formula holds:
\begin{equation}\label{formula induttiva per ptildetildeA}
\ptildetilde_A=\frac{e^{g_1\ptildetildeAtras}-1}{g_1}
\end{equation}
where $\ptildetildeAtras$ is $\ptildetilde_A(g_1,g_2,g_3,\cdots)$ evaluated in $(g_2,g_3,g_4,\cdots)$. 
\end{teo}
\proof
The formula is recursive, as one can easily check   by induction on the valency.   

Let $\T$ be a tree on $n\geq 2$ leaves (recall that we are taking into account only minimal
representatives modulo trivial tails). Let $i_1^{m_1}\cdots i_r^{m_r}$ be a partition of $n$ of length $k$ made by
positive integers $i_1,\cdots,i_r$ such that, for each $j\in\{1,\cdots,r\}$, $i_j$ occurs $m_j$ times (and
$k=\sum_{j=1}^r m_j$). Suppose that  exactly $k$ edges stem from the root; furthermore, suppose   that, if we cut off the root of $\T$ and these   edges,  we get a forest of $k$ trees,
$\{\T_1,\cdots,\T_k\}$, such that, for each $j\in\{1,\cdots,r\}$, $m_j$ of them are isomorphic and have   $i_j$ leaves (here
we are considering also the degenerate tree given by a single leaf).\\
If, for every $i\in\{1,\cdots,k\}$, we call $m_{\T_i}$ the monomial of $\T_i$ in $\ptildetilde_A$ and $m_\T$ the one of
$\T$ we have:
$$m_\T=g_1^{k-1}\frac{1}{m_1!m_2!\cdots m_r!}\prod_{i=1}^k m_{\T_i}\left[1\right].$$
We conclude by observing that $\prod_{i=1}^k m_{\T_i}\left[1\right]$ appears exactly $\frac{k!}{m_1!m_2!\cdots m_r!}$ times
in $\left(\ptildetildeAtras\right)^k$.\\ 
\qed

\subsection{Some Examples}
\label{sec:examples}

Theorem \ref{formula induttiva per ptildetildeA teorema} allows us to compute $\ptildetilde_A$; once we have
$\ptildetilde_A$ we can compute $\ptilde_A$ and, using theorem \ref{serie da ptilde}, the series $\phi_A$; here we
exhibit some examples of these computations made with the help of the Computer Algebra system Axiom.\\
\\
As a first example, we show the monomials  of $\ptildetilde_A$ of valency less than or equal to 3 and degree less than or equal to 3:
{\small \begin{align*}
\ensuremath{& \frac{1}{24} \  {g_3^3}+\left(  \frac{7}{24} \  g_2+
\frac{7}{24} \  g_1+ \frac{1}{6}
\right)\  {g_3^2}
+\left(  \frac{1}{4} \  {g_2^2}+\left(  \frac{3}{4} \ g_1+ \frac{1}{2}\right)
\  g_2+ \frac{1}{4} \  {g_1^2}+ \frac{1}{2} \  g_1+ \frac{1}{2}\right)\  g_3+ \\}
\ensuremath{& +\frac{1}{24} \  {g_2^3}+\left( \frac{7}{24} \  g_1+\frac{1}{6}\right)
\  {g_2^2}+\left( \frac{1}{4} \  {g_1^2}+ \frac{1}{2} \  g_1+\frac{1}{2}\right)
\  g_2+ \frac{1}{24} \  {g_1^3}+\frac{1}{6} \  {g_1^2}+
\frac{1}{2} \  g_1+1. }
\end{align*}}
If, for example, we look at terms of degree 3 we have
$$\frac{g_3^3}{4!}+\frac{g_2g_3^2}{8}+\frac{g_2g_3^2}{6}+\frac{g_2^2g_3}{4}+\frac{g_1g_3^2}{8}+\frac{g_1g_3^2}{6}+\frac{
g_1g_2g_3}{2}+\frac {g_1g_2g_3}{4}+\frac{g_1^2g_3}{4}.$$
The nine  monomials correspond  to the  levelled trees on 4 leaves with 3 levels (modulo equivalence):
$$
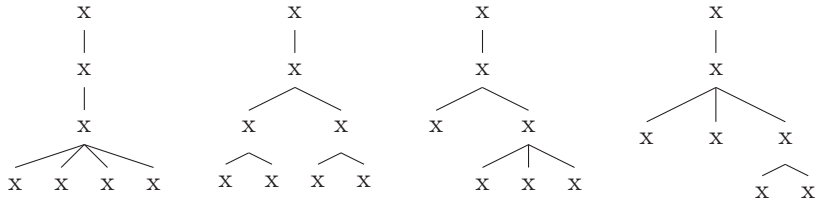
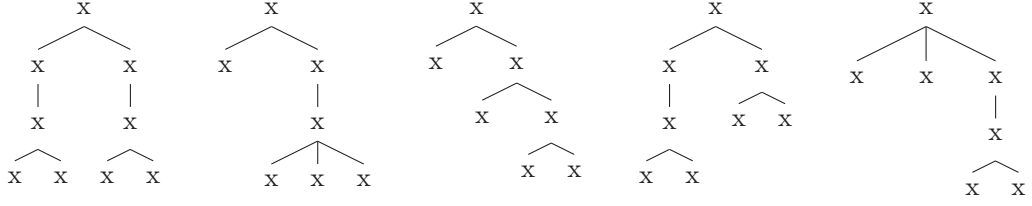
\qtreecenterfalse
\Tree[ .x [ [ x x x x ].x ].x ]
\hskip 0.3in
\Tree[ .x [ [ x x ].x [ x x ].x ].x ]
\hskip 0.3in
\Tree[ .x [ x [ x x x ].x ].x ]
\hskip 0.3in
\Tree[ .x [ x x [ x x ].x ].x ]
$$
$$
\qtreecenterfalse
\Tree[ .x [ [ x x ].x ].x [ [ x x ].x ].x ]
\hskip 0.3in
\Tree[ .x x [ [ x x x ].x ].x ]
\hskip 0.3in
\Tree[ .x x [ x [ x x ].x ].x ]
\hskip 0.3in
\Tree[ .x [ [ x x ].x ].x [ x x ].x ]
\hskip 0.3in
\Tree[ .x x x [ [ x x ].x ].x ]
$$
\bigskip

Then  we show $\ptilde$  up to degree 5, without bad monomials: 
\begin{align*}
\intertext{Degrees 1, 2 and 3}
\ensuremath{& g_0 + { \frac{1}{2}} \  g_0 \  g_1}+{{ \frac{1}{2}} \  {g_0^2}}+  {{ \frac{1}{2}} \  g_0 \  g_1 \  g_2}+{{ \frac{1}{6}} \  g_0 \  {g_1^2}}+{{
\frac{1}{2}} \  {g_0^2} \  g_1}+{{ \frac{1}{6}} \  {g_0^3}}
\intertext{Degree 4}
\ensuremath{&  \frac{3}{4} \  g_0 \  g_1 \  g_2 \  g_3+ \frac{7}{24} \  g_0 \  g_1 \
{g_2^2}+\left( \frac{1}{4} \  g_0 \  {g_1^2}+ \frac{3}{4} \  {g_0^2} \ g_1 \right)
\  g_2+ \frac{1}{24} \  g_0 \  {g_1^3}+ \\}
\ensuremath{&
 +\frac{7}{24} \  {g_0^2} \
{g_1^2}+ \frac{1}{4} \  {g_0^3} \  g_1+ \frac{1}{24} \  {g_0^4} }
\intertext{Degree 5}
\ensuremath{&  \frac{3}{2} \  g_0 \  g_1 \  g_2 \  g_3 \  g_4+ \frac{5}{8}\  g_0 \  g_1 \
g_2 \  {g_3^2}+\left(  \frac{7}{12} \  g_0 \  g_1 \  {g_2^2}+\left(
\frac{1}{2} \  g_0 \  {g_1^2}+ \frac{3}{2} \  {g_0^2} \  g_1 \right)\  g_2 \right)
\  g_3+ \\}
\ensuremath{& +\frac{1}{8} \  g_0 \  g_1 \  {g_2^3}+\left(  \frac{5}{24} \
g_0 \  {g_1^2}+ \frac{5}{8} \  {g_0^2} \  g_1 \right)
\  {g_2^2}+\left(  \frac{1}{12} \  g_0 \  {g_1^3}+ \frac{7}{12} \
{g_0^2} \  {g_1^2}+ \frac{1}{2} \  {g_0^3} \  g_1 \right)
\  g_2+  \\}
\ensuremath{& +\frac{1}{120} \  g_0 \  {g_1^4}+ \frac{1}{8} \  {g_0^2} \
{g_1^3}+ \frac{5}{24} \  {g_0^3} \  {g_1^2}+ \frac{1}{12} \  {g_0^4} \  g_1+ \frac{1}{120} \  {g_0^5} }
\end{align*}
At last here it is the series $\phi_A$ up to degree 7 (with respect to $t$): 
\begin{align*}
\ensuremath{& \phi_A(q,t)=t+ \frac{1}{2} \ t^2 + \left( \ \frac{1}{6} \ q+\frac{1}{6} \ \right) \ t^3 + \ \left( \
\frac{1}{24} \ q^2+\frac{1}{3} \ q+\frac{1}{24} \ \right) \ t^4+ \\}
\ensuremath{& +{{\left( {{ \frac{1}{{120}}} \  {q^3}}+{{ \frac{{41}}{{120}}} \  
{q^2}}+{{ \frac{{41}}{{120}}} \  q}+{ \frac{1}{{120}}} 
\right)}
\  {t^5}}+\\}
\ensuremath{&+{{\left( {{ \frac{1}{{720}}} \  {q^4}}+{{ \frac{{187}}{{720}}} \  
{q^3}}+{{ \frac{{61}}{{60}}} \  {q^2}}+{{ \frac{{187}}{{720}}} \  q}+{ 
\frac{1}{{720}}} 
\right)}
\  {t^6}}+\\}
\ensuremath{&{+\left( {{ \frac{1}{{5040}}} \  {q^5}}+{{ \frac{{19}}{{112}}} \  {q^4}}+{{ 
\frac{{2389}}{{1260}}} \  {q^3}}+{{ \frac{{2389}}{{1260}}} \  {q^2}}+{{ 
\frac{{19}}{{112}}} \  q}+{ \frac{1}{{5040}}} 
\right)}
\  {t^7}+\cdots}
\end{align*}

\subsection{Type $B_n$}

Let $V$ be a complex vector space of dimension $n$. We consider   in $V^*$  the line arrangement  corresponding to a  root system of type $B_n$ and denote by $\maxBn$
its closure under the sum and by $\F_{B_n}$ the set of irreducible subspaces  of $\maxBn$. Our aim is  to compute the series:
\begin{equation}\label{serie poincare' Bn}
 \phi_B(q,t):=\frac{t}{2}+\sum_{n\geq 2} P_\maxBn(q)\frac{t^n}{2^nn!} \ \in \ \Q[q][[t]]
\end{equation}
where, for each $n\geq 2$, $P_\maxBn(q)$ is the Poincar\'e polynomial of $Y_\maxBn$.\\
In \cite{YuzBasi}, Yuzvinsky divided the elements of $\F_{B_n}$ in two classes: strong elements and weak elements; if we
call $x_1,\cdots,x_n \ \in \ V^*$ the coordinate functions, strong elements are the subspaces of $V^*$ like
$<x_{i_1},\cdots,x_{i_k}>$ ($k\geq 1$),  whose annihilator in $V$ is the subspace $H_{i_1,\cdots,i_k}:=\{x_{i_1}=\cdots
=x_{i_k}=0\}$. They can be put in bijective correspondence with subsets of $\{1,\cdots,n\}$ of cardinality greater then
or equal to $1$ (such subsets will be called strong). A weak element is a subspace  whose annihilator is of type $L_{i_1,\cdots,i_k,j_1,\ldots, j_s}:=\{x_{i_1}=\cdots
=x_{i_k}=-x_{j_1}=\cdots =-x_{j_s}\}$ (\(r+s\geq 2\)); therefore weak elements  can be put in a bijective correspondence with subsets of $\{1,\cdots,n\}$ of cardinality greater than or equal
to $2$ equipped with a partition (possibly trivial)  into  $2$ parts  (such subsets will be called weak).\\
Moreover,  if we order $\F_{B_n}$ by inclusion of subspaces, we can read  this order as
follows:
\begin{itemize}
 \item a  subset  that includes a strong subset of  $\{1,\cdots,n\}$ is strong;
 \item a weak subset $A$ is smaller than a strong subset $B$ if and only if  $A\subset B$;
 \item a weak subset $A=A_1\cup A_2$ is smaller than a weak subset $B=B_1\cup B_2$ if  and only if  either
$A_i\subset B_i$ ($i=1,2$) or $A_1\subset B_2$ and $A_2\subset B_1$.
\end{itemize}
Coming to the maximal building set, we observe that there is  a bijective correspondence between elements of
$\maxBn$ and families of disjoint subsets of $\{1,\cdots,n\}$ in which at most one is  strong and in each of the
weak ones a partition into two parts  is fixed. Given two such families $X=\{X_1,\cdots,X_k\}$ and
$Y:=\{Y_1,\cdots,Y_h\}$ we say that $X$ is greater than $Y$ (and write $X\supset Y$) if for every $i\in\{1,\cdots,h\}$
there exists $j\in\{1,\cdots,k\}$ such that $Y_i\subset X_j$ as elements of $\F_{B_n}$.

It is again possible to associate levelled forests to $\maxBn$-nested sets (a $\maxBn$-nested  is a   subset of $\maxBn$ strictly ordered by inclusion);  the rules are the same as in the case $A_{n-1}$ but
we have to divide the vertices of our graphs into  two classes: weak vertices and strong vertices. We notice that  we lose
the information concerning partitions of weak sets. From now we call "strong tree" a tree with at least one strong
vertex and "weak tree" a tree with no strong vertices;   a forest is ``strong'' if it contains a strong tree, otherwise is weak.

Let now $\Sc$ be a $\maxBn$-nested set and $\Gamma\left(\Sc\right)$ be its associated forest; 
 let us denote by \(B\) (resp. \(A\)) the element of \(\Sc\) determined by the vertices (not leaves) at level $k$ (resp. $k+1$); then   $A$ is the maximal element of $\Sc$ strictly contained
in $B$.  If \(B\)  is given by a family of  weak subsets, \(d_{\{A\},B}\) can be  computed, in terms of outgoing edges,  exactly  as in the \(A_n\) case.
Otherwise,  \(B\) is associated to a family  $\{B_1,\cdots,B_k\}$ (\(k\geq 1\)) of subsets of $\{1,\cdots,n\}$, where  \(B_1\) is  strong.
Then we have 
\begin{equation}\label{bound coomologia Bn}
 d_{\{A\},B}=|\wout(v_{B_1})|+\sum_{i=2}^{k} \left(|\out(v_{B_i})|-1\right)
\end{equation}
where $v_{B_i}$ is the vertex of $\Gamma\left(\Sc\right)$ which corresponds to  $B_i$ and $\wout(v_{B_1})$ is the set of outgoing
edges from $v_{B_1}$ to  a weak vertex (we are considering the leaves as weak vertices).

The following lemma and corollary explain how to take in account the information on the partitions
associated to weak sets, which is not contained in the graphs.

\begin{lemma}\label{numero di nested che corrispondono a una foresta nel caso Bn irriducibile}
 (see \cite{YuzBasi}).\\
Let $\Sc$ be  a $\F_{B_n}$-nested set and $\Gamma\left(\Sc\right)$ be its  associated
forest.  If we denote by $\pi\left(\Gamma\left(\Sc\right)\right)$ the number of different
$\F_{B_n}$-nested sets  $\U$ such that $\Gamma\left(\U\right)=\Gamma\left(\Sc\right)$,  then
$$\log_2 \pi\left(\Gamma\left(\Sc\right)\right)=\sum_{v_B} \dim B$$
where \(v_B\) ranges over  all the maximal  weak vertices (not leaves), i.e. the weak vertices (not leaves) which are not preceded, according to the orientation,  by other weak vertices) .
\end{lemma}
\begin{cor}\label{numero di nested che corrispondono a una foresta etichettata}
 Let $\Gamma=\Gamma\left(\Sc\right)$ be a levelled  forest associated to a $\maxBn$-nested set $\Sc$. Let 
$\{v_{X_1},\cdots,v_{X_j}\}$ be  the maximal   weak vertices  of \(\Gamma\left(\Sc\right)\).   Then $\Gamma$
corresponds
to $2^{\sum_{i=1}^j (|X_i|-1)}$ different $\maxBn$-nested sets.
\end{cor}

As in the \(A_n\) case, to compute   the Poincar\'e series  (\ref{serie poincare' Bn}) we define a series in infinite variables \(g_0,g_1,g_2,\cdots\).  We need to extend the  definition of trivial tail to  strong trees.
\begin{defin}
 A  \textbf{trivial tail} of a levelled   oriented  strong tree $\T$  is  given by a weak subtree $\T'$
on a single leaf which stems from  a  vertex \(v\)  of $\T$  with $|\out(v)|=1$. 
\end{defin}

Then we  define a series which will take into account the contribution of strong forests to  the Poincar\'e series: 
{\small $$\Qtilde_B(g_0,g_1,g_2,\cdots):=\sum_\Gamma\frac{1}{|\Aut \left(\Gamma\right)|}\left(\frac{g_0}{2}\right)^{
\tr\left(\Gamma\right)-1}\prod_{v\in\Gamma_s}\left(\frac{g_{l(v)}}{2}\right)^{|\wout(v)|}\prod_{v\in\Gamma_w}g_{l(v)}^{
|\out(v)|-1}$$}
where $\Gamma$ runs among minimal representatives (modulo trivial tails) of  levelled oriented strong forests on $n\geq 1$ leaves, $\Gamma_s$ is the
set of strong vertices of $\Gamma$, $\Gamma_w$ is the set of  weak vertices  (not leaves); we  notice that an automorphism  sends strong vertices to  strong vertices and weak vertices  to  weak
vertices.\\
Then we put: 
\begin{equation}\label{definizione di Q}
Q_B(g_0,g_1,g_2,\cdots):=\ptilde_A(\frac{g_0}{2},g_1,g_2,\cdots)+\Qtilde_B(g_0,g_1,g_2,\cdots).
\end{equation}
\begin{rem} We have  that 
\(\ptilde_A(\frac{g_0}{2},g_1,g_2,\cdots)=e^{\frac{g_0}{2}\ptildetilde_A(g_1,g_2,\cdots)}-1\).  
\end{rem}
We  notice that, if  $\Gamma$ is a (strong or weak) levelled oriented   forest on $n\geq 2$ leaves   then the corresponding  monomial  $m_\Gamma$ in $Q_B(g_0,g_1,g_2,\cdots)$ has degree  $n$.
\begin{teo}\label{serie da Q}
If we remove bad monomials from $Q_B$ and replace $g_0$ with $t$ and, for every $i,r\geq 1$, $g_i^r$ with
$\frac{q^r-q}{q-1}t^r$, we obtain the series (\ref{serie poincare' Bn}).
\end{teo}
\proof
It is a   computation very similar to the one  of Theorem \ref{serie da ptilde}: \(\ptilde_A(\frac{g_0}{2},g_1,g_2,\cdots)\) counts the contribution of weak  forests,   \(\Qtilde_B\) of strong forests.
\qed
We now need to compute $\Qtilde_B$; to this end we define
$$\Qtildetilde_B(g_1,g_2,\cdots):=\sum_\T
\frac{1}{|\Aut\left(\T\right)|}\prod_{v\in\T_s}\left(\frac{g_{l(v)}}{2}\right)^{|\wout(v)|}\prod_{v\in\T_w}g_{l(v)}^{|
\out(v)|-1}$$
where $\T$ runs among minimal representatives of classes of strong trees on $n\geq 1$ leaves, $\T_s$ is the set of strong vertices
of $\T$ and $\T_w$ is the set of the weak ones (not leaves). Then, since each strong forest has exactly   one strong tree
we have:  $$\Qtilde_B=\left(\ptilde_A(\frac{g_0}{2},g_1,g_2,\cdots)+1\right)\Qtildetilde_B.$$

\begin{teo}\label{formula induttiva per Qtildetilde teorema}
 The following inductive formula holds:
 \begin{equation}\label{formula compatta per Qtildetilde}
\Qtildetilde_B=\ptilde_A(\frac{g_0}{2},g_1,g_2,\cdots)\left[1\right] \ + \ \Qtildetildetras \ \left(1 \ + \ \ptilde_A(\frac{g_0}{2},g_1,g_2,\cdots)\left[1\right]\right).
\end{equation}
\end{teo}
\proof
We will prove the equivalent formula
{\small \begin{equation}\label{formula induttiva per Qtildetilde}
\Qtildetilde_B=\left(\sum_{j\geq0}\frac{g_1^j\left(\frac{1}{2}\ptildetilde_A(g_1,g_2,\cdots)[1]\right)^j}{j!}\right)\Qtildetildetras+\sum_{ j\geq
1}\frac{g_1^j\left(\frac{1}{2}\ptildetilde_A(g_1,g_2,\cdots)[1]\right)^j}{j!}.
\end{equation}}
Let $\T$ be a strong tree on $n\geq 1$ leaves; suppose that $\T$ has only one strong vertex (therefore its root is strong). 
Let $i_1^{m_1}\cdots i_r^{m_r}$ be a partition of $n$ of length $k\geq 1$ made by positive integers
$i_1,\cdots,i_r$ such that, for each $j\in\{1,\cdots,r\}$, $i_j$ occurs $m_j$ times. 
Let us  suppose that $k$ edges  stem from the (strong)
root of $\T$ and   that, if we cut off the root of $\T$ and these   edges,  we get a forest of $k$ weak trees,
$\{\T_1,\cdots,\T_k\}$, such that, for each $j\in\{1,\cdots,r\}$, $m_j$ of them are isomorphic and have   $i_j$ leaves. 
If we call $m_\T$ the monomial of $\T$ in $\Qtildetilde_B$ and, for every $i\in\{1,\cdots,k\}$, $m_{T_i}$ the one of
$\T_i$ in \(\ptildetilde_A(g_1,g_2,\cdots)\), we have
$$m_\T=\frac{1}{m_1!m_2!\cdots m_r!}\frac{g_1^k}{2^k}\prod_{i=1}^k m_{\T_i}\left[1\right].$$
To obtain the second addendum on the right side of formula (\ref{formula induttiva per Qtildetilde})  it's enough to notice  that $\prod_{i=1}^k m_{\T_i}\left[1\right]$ appears exactly $\frac{k!}{m_1!m_2!\cdots m_r!}$
times in $\left(\ptildetilde_A(g_1,g_2,\cdots)[1]\right)^k$.\\
Suppose now that $\T$ has more than one strong vertex and  that the (strong) root of $\T$ is connected to $j+1\geq 1$ vertices such that $j$ of them are weak and one (which we will denote by \(v\))  is strong.\\
If $j=0$ then exactly one edge stems  from the root of $\T$  and (by assumption) it reaches  $v$. If we
call $\T'$ the tree which stems from $v$ we have that $m_\T=m_{\T'}\left[1\right]$.\\ 
Let $j>0$ and   suppose that the $j$ subtrees whose roots are  the $j$ weak vertices
 are divided into \(h\) subsets containing respectively $m_1,m_2, \ldots, m_h$  isomorphic trees.  Let us denote by 
$\{\T_1,\cdots,\T_{j}\}$ these trees and let   $\T_{j+1}$ be the  strong tree whose root is \(v\). We have
$$m_\T=\frac{1}{m_1!m_2!\cdots m_h!}\frac{g_1^j}{2^j}m_{\T_{j+1}}\left[1\right]\prod_{q=1}^j m_{\T_q}\left[1\right]$$
where, as usual, $m_\T$ is the monomial of $\T$ in $\Qtildetilde_B$, $m_{\T_{j+1}}$ is the one of $\T_{j+1}$ and, for
each
$q\in\{1,\cdots,h\}$, $m_{\T_q}$ is the monomial of $\T_q$ in \(\ptildetilde_A\).\\
We end by observing that $\prod_{i=1}^km_{\T_i}\left[1\right]$ appears exactly $\frac{j!}{m_1!m_2!\cdots m_h!}$ times in
$\ptildetilde_A^j$.\\
\qed 

\subsection{Type $D_n$}
\label{subsection:Dn}
Let $V$ be a complex vector space of dimension $n$. We consider   in $V^*$  the line arrangement  corresponding to a  root system of type  $D_n$  and denote by $\maxDn$
its closure under the sum.\\
The series we are interested in is the following:
\begin{equation}\label{serie poincare' Dn}
 \phi_D(q,t):=t+q\frac{t^2}{4}+\sum_{n\geq 3}P_{\maxDn}(q)\frac{t^n}{n!2^{n-1}} \ \in \ \Q[q][[t]]
\end{equation}
where, as usual, for each $n\geq 3$, $P_{\maxDn}(q)$ is the Poincar\'e polynomial of $Y_\maxDn$. 
The combinatorics is essentially the same as in the case $B_n$: the only difference is that strong sets must have
cardinality at least two and then, for the computation of the Poincar\'e series, we just need to modify a little what we
have done in that case.\\
We set
{\small $$\Qtilde_D(g_0,g_1,g_2,\cdots):=2\sum_\Gamma\frac{1}{|\Aut \left(\Gamma\right)|}\left(\frac{g_0}{2}\right)^{
\tr\left(\Gamma\right)-1}\prod_{v\in\Gamma_s}\left(\frac{g_{l(v)}}{2}\right)^{|\wout(v)|}\prod_{v\in\Gamma_w}g_{l(v)}^{
|\out(v)|-1}$$}
where  $\Gamma$ is a  strong  levelled oriented forest whose strong vertices correspond to subsets of cardinality   least two and
$$\Qtildetilde_D(g_1,g_2,\cdots):=2\sum_\T
\frac{1}{|\Aut \left(\T\right)|}\prod_{v\in\T_s}\left(\frac{g_{l(v)}}{2}\right)^{|\wout(v)|}\prod_{v\in\T_w}g_{l(v)}^{|
\out(v)|-1}$$
where the strong  vertices of the strong  levelled oriented tree $\T$ correspond to subsets of cardinality   least two.
Now, if we define 
$Q_D:=2\ptilde_A(\frac{g_0}{2},g_1,g_2,\cdots)+\Qtilde_D$ we can compute the Poincar\'e series in the same
way described in theorem \ref{serie da Q}; moreover we have that 
$$\Qtilde_D=\left(\ptilde_A(\frac{g_0}{2},g_1,g_2,\cdots)+1\right)\Qtildetilde_D$$
and $\Qtildetilde_D$ satisfies  a  recurrence relation similar to (\ref{formula induttiva per Qtildetilde}):
{\small \[\label{formula induttiva per QtildetildeD}
\Qtildetilde_D=\left(\sum_{j\geq0}\frac{g_1^j\left(\frac{1}{2}\ptildetilde_A(g_1,g_2,\cdots)[1]\right)^j}{j!}\right)\Qtildetilde_D[1]+2\left( \sum_{ j\geq
1}\frac{g_1^j\left(\frac{1}{2}\ptildetilde_A(g_1,g_2,\cdots)[1]\right)^j}{j!}\right) -g_1.
\]}

\section{Induced subspace arrangements}
\label{sec:final}

The tensor product allows us to obtain  new  building  subspace arrangements  $\G_h$ starting from a given building   
arrangement $\G$ in  \(V^*\).
\begin{defin}
 We will call `induced by $\G$' the 
subspace arrangement $\G_h$ in $V^*\otimes \C ^h$ (\(h\geq 1\)) given by the 
subspaces $A \otimes \C^h$, as $A$ varies in $\G$. 
\end{defin}
For instance,  if  \(\G\) is a building set associated to a root system of type \(A\),  the complements of the arrangements $\G_h$ are classical generalizations of the pure braid space  (see  \cite{cohen1976hil} and  \cite{sundaramwelker}).

It is immediate to check that, for  any  given building   
arrangement $\G$ in  \(V^*\), \(\G_h\) is still building, and therefore one can consider its   De Concini-Procesi model.
Let us focus on the case when  the starting arrangements are the  maximal building sets of type \(A,B (=C), D\).  Our series \(\ptilde_A, Q_B, Q_D\) allow us to obtain    quickly  the  Poincar\'e series of  the families of models associated to the   induced  building  sets:  we only have  to perform different substitutions for the variables \(g_0,g_1,\ldots \)

For instance, let us fix \(h\geq 1\) and consider  the \(A\) case:  after  removing  bad monomials from   \(\ptilde_A\),   if we  replace  $g_0$ with $t$ as before and, for every $i,r\geq 1$,  $g_i^r$ with
$\frac{q^{rh}-q}{q-1}t^r$, we obtain the Poincar\'e series for the models \(Y_{(\Cc_{A_n})_h}\) (the same substitutions work also in the other cases).

\addcontentsline{toc}{section}{References}
\nocite{*}
\bibliographystyle{acm}
\bibliography{Bibliogpre} 
\end{document}